\documentclass[12pt]{amsart}
\usepackage{lscape}
\usepackage{color}

\newtheorem{thm}{Theorem}[section]

\newtheorem{conj}[thm]{Conjecture}

\theoremstyle{definition}

\theoremstyle{remark}

\numberwithin{equation}{section}

\def\dd{\mathinner{\ldotp\ldotp}}
\newcommand{\PP}{\mathcal{P}}
\newtheorem{qn} [thm] {Question}

\title{Solved and Unsolved Problems about Abelian Squares}

\begin{document}

\author{Jamie Simpson}%
\address{Department of Mathematics and Statistics \newline
\indent Curtin University of Technology\newline \indent GPO Box
U1987 \newline
\indent Perth, Western Australia 6845\newline\indent Australia}%
\email{Jamie.Simpson@maths.curtin.edu.au}

\begin{abstract} We present and discuss a number of known results and open problems abelian squares in words on small alphabets.
\end{abstract}

\maketitle  We use the usual notation for combinatorics on
words. A word of $n$ elements is $w=w[1 \dd n]$, with $w[i]$ being
the $i$th element and $w[i\dd j]$ the \emph{factor} of elements from
position $i$ to position $j$.  If $i=1$ then the factor is a
\emph{prefix} and if $j=n$ then it is a \emph{suffix}. The letters
in $w$ come from some \emph{alphabet} $A$. The set of all finite
words with letters from $A$ is $A^*$. The \emph{length} of $w$,
written $|w|$, is the number of occurrences of letters in $w$ and
the number of occurrences of the letter $a$ in $w$ is $|w|_a$. If $w$ is a word whose letters come from an alphabet $a_1,a_2,\dots,a_t$ then the \emph{Parikh vector} of $w$, written  $\PP(w)$, is the vector $[|w|_{a_1},|w|_{a_2},\dots,|w|_{a_t}]$   An \emph{abelian square} is a factor $uv$ in which $\PP(u)=\PP(v)$. We say that abelian squares $u_1v_1$ and $u_2v_2$ are \emph{equivalent}  if $\PP(u_1)=\PP(u_2)$. Thus $abab$ and $abba$ are equivalent abelian squares but $abba$ and $aaaa$ are non-equivalent abelian squares.  We are interested in the maximum and minimum number of abelian squares that can occur in a word of length $n$.  This depends on the size of the alphabet we use (in what follows the alphabet will usually be binary) and how the squares are counted.  We can count the total number of abelian squares, the number of distinct abelian squares or the number of non-equivalent abelian squares.   This gives us many problems to consider.  We further expand the set of problems by considering the same set of extremal problems but applied to circular words rather than linear words. For circular words we do not consider abelian squares whose length is greater than the word length, that  is, we don't allow abelian squares to overlap themselves.  Changing from linear to circular words can simplify problems as we are avoiding messiness at the ends of words.  For example the minimum number of non-equivalent abelian squares in a binary circular word was established with a short straight-forward proof in \cite{FSP}, but the same question for linear words has not been settled and appears to be difficult. On the other hand the maximum possible number of distinct non-empty palindromes in a word of length $n$ is easily shown to be $n$ but the equivalent problem for circular words is more difficult \cite{JS}.  To summarise we are looking at the
\begin{equation}
 \left\{ \begin{array}{l}\text{maximum}\\ \text{minimum} \end{array}\right\} \left\{\begin{array}{l}\text{total number}\\ \text{number of distinct} \\ \text{number of non-equivalent} \end{array}\right\} \text{ abelian squares in a }  \left\{ \begin{array}{l}\text{linear}\\ \text{circular} \end{array}\right\}
\end{equation} word of length $n$ on an alphabet of size $t$.

This gives 12 problems with possible sub-problems depending on alphabet size.  It will be convenient to adopt the following shorthand for these problems.  The required bounds for a binary word of length $n$ are written $C_1C_2C_3(n)$ where $C_1$ is either $X$ or $M$ for maximum and minimum respectively, $C_2$ is $T$, $D$ or $N$ for total, distinct and non-equivalent and $C_3$ is either $L$ or $C$ for linear and circular.  Bounds for words on larger alphabets will be treated separately. We have
$$XC_2C_3(n) \ge MC_2C_3(n)$$ since maximums are at least as large as minimums,
$$C_1TC_3(n) \ge C_1DC_3(n) \ge C_1NC_3,$$ which follows immediately from the meanings of $T$, $D$ and $N$, and
$$C_1C_2C(n) \ge C_1C_2L(n)$$ since a circular word contains the same abelian squares as the corresponding linear word plus any abelian squares that straddle the ends of the word.\\

\begin{center} \textbf{The maximum number of abelian squares in a linear word}\\
\end{center}

 The maximum total number of abelian squares in a linear word of length $n$ is  $\lfloor n/2 \rfloor \lceil n/2 \rceil$ and is attained by $a^n$ so we have
 $$XTL(n)=\lfloor n/2 \rfloor \lceil n/2 \rceil.$$ For circular words we have
 $$XTC(n)=n\lfloor n/2 \rfloor $$ which is also attained by $a^n$.

For the other two ways of counting abelian squares we will only consider the case of binary words.  I'm not aware of any work done on larger alphabets and it's likely that increasing the alphabet size will not increase the maximum number of abelian squares. In fact Fici and Mignosi \cite{FM} make the following conjecture.
 \begin{conj} If a word of length $n$ contains $k$ many distinct abelian square factors, then there exists a binary word of length $n$ containing at least $k$ many distinct abelian square factors.
 \end{conj}

   Tables 1 and 2 give  values of $XDL(n)$ and $XNL(n)$ for small values of $n$.

\begin{table}[h]
 \centering
 \begin{tabular}{lll}
$n$ & \emph{XDL(n)} & Example attaining bound\\
\hline
 1&0&a\\
 2&1&aa\\
 3&1&aaa\\
 4&2&aaaa\\
 5&3&aabba\\
 6&4&aababa\\
 7&5&aababaa\\
 8&7&aabbaabb\\
 9&9&aaaabbaaa\\
 10&11&aaaabbaaaa\\
 11&13&aaaabbaaabb\\
 12&15&aaaabbaaaabb\\
 13&17&aaaaabbaaaabb\\
 14&21&aaaaaabbaaaabb\\
 15&23&aaaababaaabbaaa\\
 16&26&aaaabaabaaababaa\\
 17&30&aaaababaaabaabaaa\\
 18&34&aababaaabaabaaaabb\\
 19&38&aababaaabaabaaaabba\\
 20&43&aababaaabaabaaaabbaa\\
  \end{tabular}

\caption{\label{table1}Maximum numbers of distinct abelian squares in binary words (\emph{XDL(n)}).  The second column is sequence A262249 in the Online Encyclopedia of Integer Sequences.}
  \end{table}

\begin{table}[h]
 \centering
 \begin{tabular}{lll}
$n$& \emph{XNL(n)} & Example attaining bound\\
\hline
 1&0&a\\
 2&1&aa\\
 3&1&aaa\\
 4&2&aaaa\\
 5&3&aabbaa\\
 6&4&aabbaa\\
 7&4&aaaabba\\
 8&6&aabbaabb\\
 9&6&aaaabbaaa\\
 10&7&aaaabbaaaa\\
 11&8&aaaabbaaabb\\
 12&10&aaaabbaaaabb\\
 13&10&aaaaabbaaaabb\\
 14&11&aaaaaabbaaaabb\\
 15&12&aaaaaabbaaaaabb\\
 16&15&aaaabbbabaaabbbb\\
 17&16&aaaabbbabaaabbbba\\
 18&17&aaaabbbabaaabbbbaa\\
 19&17&aaaaaabbbabaaabbbba\\
 20&19&aaaaaabbbaabaaaabbbb\\
  \end{tabular}

\caption{\label{table1}Maximum numbers of nonequivalent abelian squares in binary words (\emph{XNL(n)}).  The second column is sequence A262265 in the Online Encyclopedia of Integer Sequences.}
  \end{table}

It is clear that both $XDL(n)$  and $XNL(n)$ are bounded above by $n^2$ since there cannot be more abelian squares in a word than there are factors. In \cite{KRRW}  Kociuska et al. show that $XDL(n) =\Theta(n^2)$ by considering the word $a^kba^kba^{2k}$, which is easily shown to contain a quadratic number of distinct abelian squares.  They mention that Gabriel Fici had previously obtained the same result by a different method, and communicated this information to them.  The number of abelian squares in this word gives a lower bound on $XDL(4k+2)$ of $$k+(\lfloor \frac{k}{2} \rfloor +1)^2+(k+1)\lceil \frac{k}{2}\rceil.$$  Although this is quadratic in $n=4k+2$ it is a fair bit smaller than the values in Table 1, for example with $k=4$ we have $n=18$ and the formula above gives a lower bound on $XDL(22)$ of 23 compared with the actual value of 34.

In the case of non-equivalent abelian squares the Kucherov et al. \cite{KRRW} show that
\begin{equation}
XDL(n)=\Omega(n^{1.5n}/\log{n})
\end{equation}
by considering words $w_k$ defined as follows:
\begin{eqnarray*}
w_1&=&ab\\
w_k&=&w_{k-1}a^kb^k \text{ for } k>1.
\end{eqnarray*} In this case the proof is distinctly non-trivial\footnote{Professor Rytter has informed me that this lower bound has been improved to $\Omega(n^{1.5})$ using a different set of words.}. The number of abelian squares in $w_k$ is difficult to estimate.  As with the distinct abelian squares case the lower bound obtained in this way is not sharp.  For example $w_4$ has length 20 and contains 13 non-equivalent abelian squares, compared with the actual value $XNL(20)=19$.  The authors of \cite{KRRW} make the following conjecture:
\begin{conj}
$$XNL(n)=O(n^{1.5}).$$
\end{conj}

In a slightly different direction Fici and Mignosi \cite{FM} consider infinite words in which the number of distinct abelian square factors of length $n$ grows quadratically with $n$ though they replace this condition with the following slightly weaker one. Write $p_w(n)$ for the factor complexity of $w$, $a(w)$ the set of abelian squares in $w$ and $F_w(n)$ the set of length $n$ factors of $w$, so that $p_w(n)=|F_w(n)|$. They say an infinite word $w$  is \emph{uniformly abelian-square rich} if there exists a positive constant $C$ such that $$\inf_{v \in F_w(n)}|a(v)| \ge Cn^2.$$  They show that the Thue-Morse word and certain Sturmian words are uniformly abelian-square rich.

\bigskip
\begin{center} \textbf{The minimum numbers of distinct and non-equivalent abelian squares in  words}\\
\end{center}
\bigskip
The situation here parallels that for non-abelian squares. \begin{table}[h]
 \centering
 \begin{tabular}{lll}
Alphabet size & Non-abelian  & Abelian \\
\hline
2 &   $\le$ 3 distinct squares & Conjectured  $\lfloor n/4 \rfloor$ distinct/non-equivalent \\
3 &  0   & Conjectured $\le 3$ distinct/non-equivalent\\
4 &  0 & 0\\
  \end{tabular}
\caption{\label{table1}Minimum numbers of distinct/non-equivalent non-abelian/abelian squares in a word of length $n$ for different alphabet sizes.}
  \end{table}

For abelian squares the case of an alphabet of size 4 is probably the most well-known.  In 1957 Erd\H{o}s \cite{E61} conjectured that a word containing no abelian square could be constructed on an alphabet of size 4.  Two authors constructed infinite abelian-square-free words on larger alphabets (size 25 by Evdokimov \cite{Ev}, size 5 by Pleasants \cite{P70}) then Ker\"{a}nen \cite{K} used an 85-uniform morphism on a four letter alphabet to construct such an infinite word.  As far as I know nobody has used a shorter morphism to produce such words. Currie and Fitzpatrick  \cite{FS} showed that there exist binary non-abelian cube-free circular words of every length.

 For a three letter alphabet M\"{a}kel\"{a} \cite{M} has made the following conjecture.
 \begin{conj}There exist infinite ternary words containing no abelian squares of length greater than two (which implies there are at most three distinct abelian squares).
  \end{conj}
  This would be analogous to the non-abelian case where it's know  that there exist infinite binary words containing only the  squares $aa$, $bb$ and $abab$. This was first proved by Fraenkel and Simpson \cite{FS95} in 1995, with much nicer proofs coming later from  Rampersad, Shallit and Wang \cite{RSW05} in 2005, Harju and Nowotka \cite{HN} in 2006 and Badkobeh and Crochemore \cite{BC} in 2011.  The longest binary words containing only zero, one and two squares have lengths, respectively, 3 (for example $aba$), 7 ($aaabaaa$) and 18 ($abaabbaaabbbaabbab$).
  M\"{a}kel\"{a}'s conjecture has been supported by  Rampersad, see \cite{FS}, who has produced a ternary word of length 3160 wherein the only abelian squares are $aa$, $bb$ and $cc$ (each occurring many times).  In the appendix there is a word of 2034 letters containing no abelian squares other than 00, 11 and 22.  This was kindly provided by Narad Rampersad.  The longest ternary words containing 0, 1 and 2 distinct abelian squares have lengths 7, 18 and 63 respectively. Examples of words attaining these bounds are
  \begin{eqnarray*}
  &&abacaba\\
  &&abcbabccacccbabcba\\
  &&abbbcbbaccbcccaccbabbbcccabbbacabacccabbbcccacbbabbbcbbaccbccca.
  \end{eqnarray*}
  The second and third examples are unique up to permutation of the alphabet.

 For a binary alphabet we must distinguish between distinct and non-equivalent abelian squares.  Recall that that $MDL(n)$ is the minimum number of distinct abelian squares in a linear binary word of length $n$ and $MNL(n)$ is the minimum number of non-equivalent abelian squares.

 For the non-equivalent case Fraenkel, Paterson and Simpson \cite{FSP} made the following conjecture.
  \begin{conj} For all positive integers $n$ $MNL(n)=\lfloor n/4 \rfloor$ and, when $n=4k+3$ the bound is only obtained by the words $a^{2k+1}ba^{2k+1}$ and $(ab)^{2k+1}a$  and their complements.  If $n$ is not congruent to 3 modulo 4 just remove 1, 2 or 3 letters from an end of the word.
\end{conj}

 In \cite{FSP} Fraenkel, Paterson and Simpson proved that a circular binary word of length $n=2k+2$ contains at least $k$ non-equivalent abelian squares and this bound is attained only be $(ab)^{k+1}$, its complement and their conjugates. Thus,

 $$MNC(2k+2)=k.$$

  In \cite{FS} Fici and Saarela  made the following conjecture.
  \begin{conj}$MDL(n)= \lfloor n/4 \rfloor$ and the only such words of length $4k+3$ containing only $k$ distinct abelian squares are $a^{2k+1}ba^{2k+1}$ and its complement.
  \end{conj}

  Since $MDL(n) \ge MNL(n)$ a proof of the Fraenkel Paterson Simpson conjecture would imply the truth of the Fici Saarela conjecture.

 In the case of distinct abelian squares in circular binary words we have the following conjecture suggested by computer experiments.

  \begin{conj} The minimum number of distinct abelian squares in a circular word of length $n$ is: \\
  (a) $(n-1)/2$ if $n$ is odd and this bound is attained only by $a^n$, $a^{n-1}b$ and their complements and conjugates.\\
  (b) $(n-2)/2$ if $n$ is even and this bound is attained only by $a^kb^{n-k}$ and its complement and their conjugates, where $k \in \{1,3,5,\dots,n-1\}$.
  \end{conj}
\bigskip
\begin{center} \textbf{The minimum total number  of abelian squares in  binary words}\\
\end{center}
\bigskip
Now we are counting the abelian squares in a word according to multiplicity, so that $ababbaaabaa$ contains a total of 7 abelian square occurrences: three copies of $aa$ and one each of $bb$, $abab$, $abba$ and $baaaba$. This is the minimum number for binary words of length 11, so $NTL(11)=7$. I'm not aware of anybody studying the following question.

\begin{qn} What is the minimum number of occurrences of abelian squares in a binary word of length $n$?
\end{qn}
  The non-abelian case has been studied by Kucherov, Ochem and Rao \cite{KOR}.  They set $m(n)$ to be the minimum number of square occurrences in a word of length $n$.   They showed that the sequence $m(n)/n$ is convergent and called the limit $\mathcal{M}$.  They exhibited an infinite word that had a density of square occurrences of $103/187$, thereby showing that $\mathcal{M} \le 103/187=0.550802\dots$ They further showed, using ideas from \cite{T02}, that $\mathcal{M} > 0.5508$.  It is remarkable that their upper and lower bounds are so close.

The following table shows the minimum number of abelian square occurrences in binary words of length $n$ for low values of $n$.

The questions discussed in this essay can be generalised in several ways and such questions have a considerable literature.  Most obviously we could consider other powers than 2. For example Dekking \cite{D} showed that you can avoid abelian cubes  with a ternary alphabet and abelian fourth powers with a binary alphabet. or abelian fractional powers.  Another variation on the theme uses the notion of $k$-equivalence introduced by Karhum\"{a}ki.  Let $|x|_u$ be the number of times the factor $u$ occurs in the word $x$.  Two words $x$ and $y$ are \emph{k-equivalent} if $|x|_u=|y|_u$ for all factors $u$ of length at most $k$.  We now say that $xy$ is a \emph{k-abelian square} if $x$ and $y$ are \emph{k-equivalent}.  An ordinary abelian square is therefore a $k$-abelian square, and  a word being a $k$-abelian square is, in general, a stronger condition than it being and abelian square and a weaker condition than it being an ordinary square. One defines $k$-abalian powers in analagously. One can avoid 3-abelian squares with a ternary alphabet and 2-abelian cubes with a binary alphabet \cite{}

\begin{table}[h]
 \centering
 \begin{tabular}{llll}
n& \emph{MTL(n)} & Example attaining & Number of binary words \\
& & bound &attaining bound \\
\hline
 1  &  0 &  a  &  2    \\
 2  &  0 &  ab  &  2    \\
 3  &  0 &  aba  &  2    \\
 4  &  1 &  abab  &  6    \\
 5  &  2 &  abbba  &  18    \\
 6  &  3 &  abbbab  &  30   \\
 7  &  4 &  abbbabb  &  52    \\
 8  &  4 &  abbbaaba  &  8    \\
 9  &  5 &  abbbaabab  &  16    \\
 10  &  6 &  abbbaaabaa  &  20    \\
 11  &  7 &  ababbaaabaa  &  12    \\
 12  &  8 &  aabaaabbbabb  &  2    \\
 13  &  10 &  abbbabbaaabaa  &  48    \\
 14  &  11 &  abbbaaabbbbaba  &  22    \\
 15  &  12 &  ababbbbaaabbbaa  &  14    \\
 16  &  14 &  abbbaaaababbbaba  &  48    \\
 17  &  15 &  ababbbabaaaabbbaa  &  4    \\
 18  &  17 &  abbbaaabaabbbbbabb  &  68   \\
 19  &  18 &  abaabbbbbabaaababbb  &  20    \\
 20  &  20 &  abbbabaaababbbbbabaa  &  122    \\
  \end{tabular}

\caption{\label{table1}Minimum numbers of  abelian square occurrences in binary words.  Note that, up to complementation and reversal the words of length 12 and 17 are unique.  The second column is sequence A268084 in the Online Encyclopedia of Integer Sequences.}
  \end{table}

\pagebreak

\bibliography{}

\pagebreak.

\begin{center} \textbf{Appendix: a ternary word containing no abelian squares other than 00, 11 and 22}\\ \end{center}

\bigskip

00010002000111000221002000111022211000111002000100222000111200022211

00221222000222110111200011220010002001110001122201112110200100221112

00012202221011122210001120011121011122211000221112220221101112221100

01200022200112022211101122202221222000100212220222111002012200010002

00011100022111002000122022210111210001112200010002011100011222101112

11100011200021220112220002212220221011211000112200012220221100020001

00022110112111000112200010022201211100020001000221002220001000200011

10002220012111000111222022212011122210200010002210022200022101112220

22212220001002221222021110112000111022210020001002211002100011121011

12201002001122211101220002110001220002220112221020011100222111220112

22022212000211122200100020012220221112000222122202221102211122211000

12022212220001002212220011121011122211100212220001002111000120222122

20012200022211000200010002220011121100022200010002000111000221011121

02211000200010222000221100011122000221220222110002211101120001002220

00111220222101112100022110001112111011122011222000102111220112220002

12220222111012022122000100020111000112200012221110212220222110112111

00011122001000222122001110211000200010002220001110002011211100011222

11101112111000112001112200010002000111220221200102211101112100012221

11220122200010022111020102211122202221222000100020111200021220011122

21102000111210111222110021220221112221100011210022212220222111222001

11022111222110022200010002011122210022200022100222122202111022201020

00111200022110212220222110111222110002001110022210111210200010002220

00111002221112011222000222120001122210022200010002000112001000222000

11102111220100020112202221002000100212021100011102111222000111022122

20002221100011121011122202221200102221110111211100200011122200122000

22211000200111222100022212220222110221112220010221101112110222001110

00221022200022111222110221110111210001122211011210002001122210020001

00221110001112101222001110002210022211122022212200010022211100012220

00100020001110002220011100011211012220010221110001120002210022212202

21110001121101112201211000222001220001000201112110002220012111

\bigskip

\end{document}